\documentclass[1p]{elsarticle}

\usepackage{hyperref}

\usepackage{amsmath,amssymb,amsfonts,amsthm,commath,mathtools}

\usepackage[T1]{fontenc}
\usepackage[utf8]{inputenc}
\usepackage[english]{babel}
\usepackage{stmaryrd}
\usepackage{booktabs}
\usepackage{array}
\usepackage{blkarray}
\usepackage{subcaption}
\usepackage{bondgraphs}
\usepackage{color}
\usepackage{graphicx}
\usepackage{natbib}
\usepackage{hyperref}
\usetikzlibrary{fit}
\usepackage{nicefrac}
\usepackage{acro}
\usepackage{pdflscape}
\usepackage{multirow}
\usepackage{comment}
\usepackage{nicefrac}

\newtheorem{definition}{Definition}
\newtheorem{property}{Property}
\newtheorem{remark}{Remark}
\newtheorem{example}{Example}
\newtheorem{theorem}{Theorem}

\definecolor{darkgreen}{rgb}{0,0.6,0.4}

\newcommand{\EE}{\mathbb{E}}
\newcommand{\FF}{\mathbb{F}}

\newcommand{\II}{\mathbb{I}}

\newcommand{\LL}{\mathbb{L}}
\newcommand{\MM}{\mathbb{M}}
\newcommand{\MMa}{\mathbb{M}_{\mathrm{a}}}
\newcommand{\NN}{\mathbb{N}}
\newcommand{\PP}{\mathbb{P}}

\newcommand{\RR}{\mathbb{R}}
\newcommand{\WW}{\mathbb{W}}

\newcommand{\cS}{\mathcal{S}}

\newcommand{\cE}{\mathcal{E}}
\newcommand{\cZ}{\mathcal{Z}}
\newcommand{\cF}{\mathcal{F}}

\definecolor{light-gray}{gray}{0.8}
\definecolor{light-green}{rgb}{0,1,0.8}

\DeclareAcronym{phs}{
  short=PHS,
  long=port-Hamiltonian system,
}

\DeclareAcronym{sp}{
  short=SP,
  long=Statistical Physics,
}

 \newcommand{\THH}[1]{\color{black}#1~\color{black}}
 
 \newcommand{\JUDY}[1]{\textcolor{black}{#1}}

\journal{Physica D}

\begin{document}

\begin{frontmatter}

\title{From Equilibrium Statistical Physics\\ Under Experimental Constraints\\ to Macroscopic Port-Hamiltonian Systems}


\author{Judy Najnudel, Thomas Hélie, David Roze, Rémy Müller}
\address{\THH{S3AM Team, Laboratory STMS (UMR 9912), IRCAM-CNRS-SU, 1 Place Igor Stravinsky, 75004 Paris, France
}}




\begin{abstract}
This paper proposes to build a bridge between microscopic descriptions of matter with internal energy, composed of many fast interacting particles inside an environment, and their port-Hamiltonian (PH) descriptions at macroscopic scale. The environment, assumed to be slow, is modeled through experimental constraints on macroscopic quantities (e.g. energy, particle number, etc), with a partitioning into two classes: non fluctuating and fluctuating values. The method to derive the PH macroscopic laws is detailed in several steps and illustrated on two standard cases (ideal gas, Ising ferromagnets). It revisits equilibrium statistical physics with a focus on this partitioning. First, the Boltzmann's principle is used to provide the statistic law of the matter. It defines a macroscopic equilibrium characterized by a scalar value, the entropy, together with thermodynamic quantities emerging from each constraint. Then, the port-Hamiltonian system is derived. The Hamiltonian (macroscopic energy) is derived as a function of the macroscopic state (entropy and the macroscopic quantities associated with the fluctuating class). The ports (flows/efforts) are related to the time-derivative of the state and the Hamiltonian gradient in a conservative way. This open system defines the reversible laws that govern standard thermodynamic quantities. Lastly, this paper presents a strategy to extend this PH system to an irreversible conservative one, given a macroscopic dissipative law.
\end{abstract}


\end{frontmatter}

\section{Introduction}

A macroscopic system (of size 10$^{-2}$ m or bigger) is constituted of matter, that is, billions of microscopic particles (of size 10$^{-9}$ m or smaller) which are collectively responsible for the system's behavior. However, studying a single particle tells nothing about the macroscopic system, just as following the trajectory of a single person is not sufficient to predict a crowd movement. Yet, solving exhaustive equations with billions of variables would be all at once much too complex and irrelevant: at a high enough scale, individual behaviors do not matter. Indeed, one is usually not interested in the particular trajectories of water molecules in one's glass, but rather in the volume, \emph{on average}, that they take. Likewise, one is not \THH{only} interested in the day's weather report, but rather in the \THH{global \emph{tendency}}.

Averages and tendencies belong to the domain of statistics, which aims to describe complex systems with a reduced number of variables. Thus, \ac{sp} computes averages on (fast) fluctuations of complex systems in order to derive (slower) macroscopic quantities, given some experimental conditions. \JUDY{Statistical arguments for the description of a system transitioning towards thermodynamic equilibrium were introduced by Ludwig Boltzmann in 1877~\cite{boltzmann1877beziehung}. This framework allows the prediction of macroscopic thermodynamic phenomena such as temperature, entropy creation, and phase transitions~\cite{landsberg2014thermodynamics}.}

\THH{Thermodynamics has been broadly studied in the context of \ac{phs}, as well as their modeling and their control (see e.g.~\cite{eberard2004port, eberard2007extension, ramirez2013irreversible, delvenne2014finite, ramirez2016passivity, van2021classical, van2021liouville}).
However,  the proper derivation of macroscopic thermodynamic PHS from complex systems with numerous degrees of freedom is seldom addressed. \JUDY{As it happens, the choice of a system representation for this kind of model reduction is all but inconsequential, and must be handled with care~\cite{gorban2006model, ottinger2007systematic}.}
In this paper, we propose a series of systematical steps in order to construct a simplified yet physically-based structured macroscopic \ac{phs} from a system that can be described by SP.
}

Note that in the scope of this work, we limit ourselves to \emph{equilibrium} \ac{sp}, in the sense that average quantities are determined for a system \emph{at thermodynamic equilibrium}, given some experimental conditions. It is compatible with studying the system dynamics, assuming that thermodynamic relaxation (the process of reaching thermodynamic equilibrium) is infinitely faster than the rate of change of experimental conditions. Based on this assumption, a macroscopic trajectory is to be understood as a succession of \THH{thermodynamic} equilibrium states.

\begin{figure}[h!]
\centering
\resizebox{\textwidth}{!}{
\begin{tikzpicture}
\node (A) [rectangle, draw, align = left] {\textbf{A. Microscopic description - Section~\ref{sec:microstate}}\\
1. $\PP$-valued particle representation\\
2. Particle configuration $\bm{m} \in \MM = \PP^\star$\\
3. Set of characterizing functions $\mathfrak{F}$};
\node (B) [rectangle, draw, align = left, below = 1cm of A] {\textbf{B. Experimental conditions - Section~\ref{sec:Ma}}\\
4. $\mathfrak{F} = \textcolor{orange}{\mathfrak{F}^\mathrm{fixed}}\cup \textcolor{purple}{\mathfrak{F}^\mathrm{free}}$\\ 
5. Fixed values $\textcolor{orange}{\bm{\theta}^\mathrm{fixed}}$\\
6. Set of accessible configurations $\MM_a(\textcolor{orange}{\bm{\theta}^\mathrm{fixed}})$};
\node (C) [rectangle, draw, align = left, below = 1cm of B] {\textbf{C. Stochastic setting - Section~\ref{sec:stoch}}\\
7. Probability distribution $p$\\ 
8. Surprisal $\cS^b_p(\bm{m})$ \THH{(in base $b$)}\\
9. Statistical entropy $\mathsf{S}^k(p)$, $k=\nicefrac{1}{\ln b}$};
\node (D) [rectangle, draw, align = left, right = 1.5cm of B]{\textbf{D. Boltzmann principle - Section~\ref{sec:proba}}\\
10. Ergodicity $\EE_p\left[\cF_i \in \textcolor{purple}{\mathfrak{F}^\mathrm{free}}\right] = \textcolor{red}{\overline \cF_i}$\\
11. Maximum entropy \THH{given}\\ 
\THH{experimental conditions \textbf{(B.)}}\\ $\rightarrow$ \THH{Thermodynamic} entropy $S^k(\textcolor{red}{\overline \cF_i})$\\
\THH{and} Lagrange mult. $\textcolor{red}{\lambda_i}$};
\node (E) [rectangle, draw, align = left, right = 1.5cm of D] {\textbf{E. Macroscopic PHS - Section~\ref{sec:final}}\\
12. \THH{Entropy to energy representation}\\$S^k(\overline \cE,\, \overline \cF_j) \leftrightarrow E(\underbrace {\overline \cS,\,\overline \cF_j}_{\bm{x}})$\\
13. Connection to ports \THH{(effort $u$, flow $y$)}\\ $\textcolor{red}{\overline \cF_j} \leftrightarrow \dot{\bm{x}}\leftrightarrow \textcolor{red}{\bm{u}}$,\\ $\textcolor{red}{\lambda_j} \leftrightarrow \nabla E(\bm{x}) \leftrightarrow \textcolor{red}{\bm{y}}$};
\draw (A.south) [->, >= latex] to (B.north);
\draw (B.east) [->, >= latex] to (D.west);
\draw (C.east) [->, >= latex] to (D.west);
\draw (D.east) [->, >= latex] to (E.west);
\end{tikzpicture}}
\caption{from equilibrium statistical physics to macroscopic port-Hamiltonian systems (PHS): method recap with the labels of the main mathematical objects introduced in each step.}
\label{fig:recap}
\end{figure}
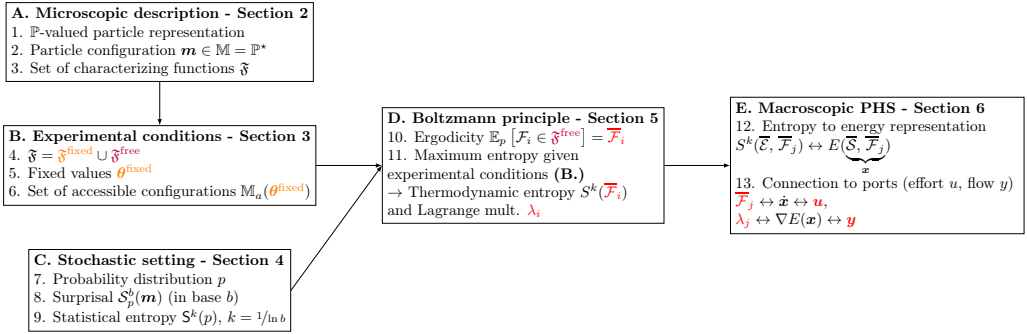

This paper is structured as follows.
In Section~\ref{sec:microstate}, we formalize the \THH{description of the microscopic configurations of} a system through the choice of (i) an ad hoc particle representation and (ii) a set of \THH{characterization functions that evaluate macroscopic quantities}.
\THH{Section~\ref{sec:Ma} addresses the experimental conditions at the macroscopic level} and their influence on the system configuration space. 
\THH{In section~\ref{sec:stoch},} we introduce a stochastic description for \THH{microscopic configurations}. 
\THH{Then, in Section~\ref{sec:proba}, we} determine the \THH{conditional} probability distribution \THH{according to the Boltzmann principle} for a system at thermodynamic equilibrium\THH{. This} allows the derivation of relevant macroscopic variables as expectations for this probability distribution.
\THH{In Section~\ref{sec:final}, we introduce the} ports \THH{and relates them to} those macroscopic variables\THH{, leading to the} macroscopic \ac{phs} model.
\THH{Finally, section~\ref{sec:FinalSummary} summarizes the practical sequence of these steps to derive the macroscopic \ac{phs} model from the microscopic description.
In addition, it presents how to derive a conservative irreversible \ac{phs} model from an additional macroscopic dissipation law.}
\\
\THH{All the steps are detailed in the following sections, as recapped in Fig.~\ref{fig:recap}.}

\section{Microstate of a system}
\label{sec:microstate}

\subsection{Particle representation \THH{($p\in\PP$)}}
In order to describe a system at a microscopic level, each of its particles must be described in a relevant way.
Depending on the system under study, one may choose to represent a particle by its position, momentum, charge, magnetic moment, etc.

\THH{
\begin{definition}[Particle set $\PP$]\label{def:P}
    Given a chosen representation to encode a particle state, we denote $\PP$ the set of all its possible values.
\end{definition}
}

\begin{example}[Particle represented by its position and momentum]
For a particle chosen to be represented by its position in space $\bm{r} \in \RR^3$ and momentum $\bm{p} \in \RR^3$, such as in a gas, \THH{the particle set is defined as} $\PP = \RR^3 \times \RR^3$. 
\end{example}

\begin{example}[Particle represented by its magnetic moment]
For a particle chosen to be represented by its magnetic moment $s \in \{-1,\, 1\}$ such as in the Ising model~\cite{ising1925beitrag, strecka2015brief}, \THH{the particle set is defined as} $\PP = \{-1,\, 1\}$. 
\end{example}

\subsection{Configuration space \THH{($\textbf{m}\in\MM$)}}
\THH{As an element of $\PP$ represents the state of one particle, a natural way to represent a configuration of particles is to concatenate elements of $\PP$.}
By analogy with formal language theory~\cite{hopcroft2001introduction}, a particular configuration of particles is \THH{chosen to be encoded as} a \emph{word} over the \emph{alphabet} $\PP$ \THH{(see remark~\ref{rk:combinatorics1} for other choices)}.

\THH{
\begin{definition}[Encoding space $\WW$]\label{def:WW}
    \color{black}
We denote $\WW := \PP^\star$ the \THH{space of encodable configurations (or, for short, the encoding space)}, where $\star$ is the Kleene operator defined by
\begin{subequations}
\begin{align}
\PP &= \{\epsilon\}, & \PP^{i + 1} &=
\left\{p_1 \boldsymbol{\cdot} p_2 \;|\; (p_1, p_2) \in \PP^i \times \PP\right\} \quad \forall i \geq 0,\\
\PP^\star &= \bigcup_{i \geq 0} \PP^i,
\end{align}
\end{subequations}
with $\epsilon$ the empty configuration and $\boldsymbol{\cdot}$ the concatenation operation.
\end{definition}
}

\begin{property}[
\THH{$\WW$} is a monoid]
By construction, the \THH{encoding space} $\WW$ is a monoid (see Def.~(\ref{def:monoid})) with associative binary operation~$\boldsymbol{\cdot}$ (concatenation) and identity element $\epsilon$ (empty configuration).
\end{property}
\begin{definition}[Monoid]
A set $\mathbb{S}$ is a monoid if it is equipped with an associative binary operation $\boldsymbol{\cdot} :\mathbb{S} \times \mathbb{S}\mapsto \mathbb{S}$ and identity element $\epsilon$, \THH{such} that for all $(s_1, s_2, s_3) \in \mathbb{S}^3$, the following properties hold
\begin{enumerate}
\item $s_1\boldsymbol{\cdot} (s_2\boldsymbol{\cdot}s_3) = (s_1\boldsymbol{\cdot} s_2)\boldsymbol{\cdot} s_3$,
\item $\epsilon\boldsymbol{\cdot}s_1 = s_1\boldsymbol{\cdot} \epsilon = s_1$.
\end{enumerate}
\label{def:monoid}
\end{definition}

\begin{figure}[tb]
\centering
\begin{subfigure}[b]{\textwidth}
\centering
\begin{tabular}{cccc}
\begin{tikzpicture}[scale = 0.6]
\draw (0, 0) -- (2, 0)
node (1) [very near start] {$\bullet$}
node [very near start, below right] {1}
node (2) [very near end] {$\bullet$}
node [very near end, below right] {2};
\draw (1)++(0, 0.7) [<-, very thick, >=latex, color=red] to (1);
\draw (2)++(0, 0.7) [<-, very thick, >=latex, color=red] to (2);
\end{tikzpicture}
&
\begin{tikzpicture}[scale = 0.6]
\draw (0, 0) -- (2, 0)
node (1) [very near start] {$\bullet$}
node [very near start, below right] {1}
node (2) [very near end] {$\bullet$}
node [very near end, below right] {2};
\draw (1)++(0, 0.7) [<-, very thick, >=latex, color=red] to (1);
\draw (2)++(0, 0.7) [->, very thick, >=latex, color=blue] to (2);
\end{tikzpicture}
&
\begin{tikzpicture}[scale = 0.6]
\draw (0, 0) -- (2, 0)
node (1) [very near start] {$\bullet$}
node [very near start, below right] {1}
node (2) [very near end] {$\bullet$}
node [very near end, below right] {2};
\draw (1)++(0, 0.7) [->, very thick, >=latex, color=blue] to (1);
\draw (2)++(0, 0.7) [->, very thick, >=latex, color=blue] to (2);
\end{tikzpicture}
&
\begin{tikzpicture}[scale = 0.6]
\draw (0, 0) -- (2, 0)
node (1) [very near start] {$\bullet$}
node [very near start, below right] {1}
node (2) [very near end] {$\bullet$}
node [very near end, below right] {2};
\draw (1)++(0, 0.7) [->, very thick, >=latex, color=blue] to (1);
\draw (2)++(0, 0.7) [<-, very thick, >=latex, color=red] to (2);
\end{tikzpicture}
\end{tabular}
\caption{Examples of microstates for a system of two particles described by their spin $s \in \{-1 \text{ (blue)}\,, 1 \text{ (red)}\}$.}
\label{fig:spin}
\end{subfigure}
\begin{subfigure}[b]{\textwidth}
\centering
\begin{tabular}{cc}
\begin{tikzpicture}[scale = 0.8]
\node [circle, draw, red] (a) at (0, 0) {};
\draw (a) [->, red] -- ++ (0.4, 0.4);
\node [circle, draw, red] (b) at (1, 1) {};
\draw (b) [->, red] -- ++ (0, -0.56);
\node [circle, draw, red] (c) at (0.2, 2) {};
\draw (c) [->, red] -- ++ (-0.4, 0.4);
\draw (-0.5, -0.5) [thick, rounded corners = 5pt] rectangle (1.5, 2.5);
\end{tikzpicture}
&
\begin{tikzpicture}[scale = 0.8]
\node [circle, draw, red] (a) at (0, 1) {};
\draw (a) [->, red] -- ++ (0, 0.56);
\node [circle, draw, red] (b) at (0.5, 0.2) {};
\draw (b) [->, red] -- ++ (-0.56, 0);
\node [circle, draw, red] (c) at (1, 2) {};
\draw (c) [->, red] -- ++ (-0.4, -0.4);
\draw (-0.5, -0.5) [thick, rounded corners = 5pt] rectangle (1.5, 2.5);
\end{tikzpicture}
\end{tabular}
\caption{Examples of microstates for a system with three  particles described by their position (circle) and momentum (arrow).}
\label{fig:tank}
\end{subfigure}
\caption{Examples of microstates for different systems.}
\end{figure}
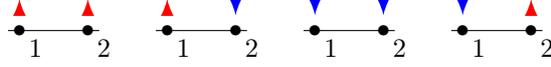
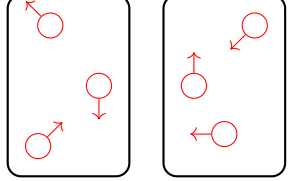

Based on the chosen representation, some configurations may not be \THH{physically} admissible~\footnote{For instance, if the chosen representation assigns a unique label to each particle, configurations in which several particles share the same label are not admissible \THH{(see remark~\ref{rk:combinatorics1}(iii) for more more details)}.}\THH{, therefore we introduce the set of microstates as follows.}
\THH{
\begin{definition}[Admissible configuration set $\MM$ and microstate $\bm{m}\in\MM$]
    We
    \color{black}
     denote $\MM \subseteq \WW$ the set of \THH{encodable configurations that are also physically admissible}.
    An element $\bm{m} \in \MM$ is called a \emph{microstate} of the system.
\end{definition}
}
In the following, 
we choose $\MM = \WW$ \THH{(see remark~\ref{rk:combinatorics1}(i) for an interpretation and (ii-iii) for examples with $\MM \neq \WW$)}.
Figure~\ref{fig:spin} shows examples of microstates for a system of particles described by their spin, and Fig.~\ref{fig:tank} shows a system of particles described by their position and momentum.

\begin{property}[\THH{$\MM$} is measurable]\label{prop:measurableM}
The pair $\left(\MM, \mathcal{P}(\MM)\right)$ where $\mathcal{P}(\MM)$ denotes the powerset of $\MM$ is a measurable space, that is, it verifies
\begin{enumerate}
\item $\MM \in \mathcal{P}(\MM)$,
\item $\mathcal{P}(\MM)$ is closed under complements: $\MM \backslash \mathbb{P} \in \mathcal{P}(\MM), \quad \forall \mathbb{P} \in \mathcal{P}(\MM)$,
\item $\mathcal{P}(\MM)$ is closed under countable unions: $\bigcup_{i = 1}^\infty \mathbb{P}_i \in \mathcal{P}(\MM) \quad \forall \mathbb{P}_1, \mathbb{P}_2, \hdots \in \mathcal{P}(\MM)$.
\end{enumerate}
\label{prop:measure}
\end{property}

\THH{
As mentioned above, the set $\MM$ could be defined on encoding spaces generated with operators other than the concatenation.
Examples are outlined in the following remark~\ref{rk:combinatorics1}.
}
\THH{
\begin{remark}[Examples of combinatorial structures and interpretations]\label{rk:combinatorics1}
The monoid $\WW$ provides a simple and natural way to encode microscopic configurations, which consists in choosing a prioritized and distinguishable representation of particles (see specification (i) below). 
But the configuration encoding can be addressed by using any appropriate combinatorics of particles, possibly choosing other specifications (see e.g.~\cite[\S\,I.2]{flajolet2009analytic} for details on operators {\sc Seq}, {\sc Mset}, {\sc Pset}, etc., mentioned below):
\begin{itemize}
    \item[(i)] \emph{Distinguishability with prioritization.}
    A word $p_i \boldsymbol{\cdot} p_{i-1} \boldsymbol{\cdot} \dots \boldsymbol{\cdot} p_1 \in \PP^i \subset \WW$ can be interpreted as describing the state values $p_{n}\in\PP$ of particles number $n=1,\dots,i$. 
    In this sense, definition~\ref{def:WW} encodes a physics with \emph{distinguishable} particles and with a \emph{priority ordering} on involved particles (particle~1 can be encoded alone, particle~2 only if 1 is involved, etc). 
    From the combinatorics point of view, the encoding space $\WW=\PP^\star$ corresponds to the sequence construction, also denoted {\sc Seq}$(\PP)$ in \cite{flajolet2009analytic}.
    \item[(ii)] \emph{Undistinguishability.}
    A natural encoding of a physics described with \emph{undistinguishable} particles is the multiset $\WW_{\mathrm{(ii)}}=${\sc Mset}$(\PP)$ composed of all the finite sets of $\PP$ (the order between elements does not count), in which arbitrary finite repetitions of elements are allowed.
    Examples of admissible configuration sets $\MM_\mathrm{(ii)}$ are $\WW_{\mathrm{(ii)}}$, or the powerset {\sc Pset}$(\PP)\subset${\sc Mset}$(\PP)$ (no repetition allowed) if the physics under consideration forbids two particles of a microstate to be excited by the same state value.
    %
    \item[(iii)] \emph{Distinguishability without priorization.}
    A means other than (i) to encode configurations with distinguishable particles is to add a distinct label to each particle such as its number $n\in\LL\subseteq\NN$.
    In this case, the encoding space can be the powerset $\WW_{\mathrm{(iii)}}=${\sc Pset}($\PP_{labeled})$, where the set of labeled particles $\PP_{labeled}$ is the cartesian product $\LL\times\PP$.
    An admissible configuration set that forbids the particle ubiquity or replication is a subset $\MM_{(iii)}$ of $\WW_{\mathrm{(iii)}}$ for which all the elements make the number of occurrences of each label $n\in\LL$ be 0 or 1.
    \item[(iv)] \emph{Complex structures.} More generally, the encoding space can be composed of complex elements, with structures involving sequences, cycles (for e.g. aromatic molecules), trees (etc.) and their combination according to precise combinatorial specifications (see \cite[I.2.3]{flajolet2009analytic}.
\end{itemize}
\end{remark}
}

\subsection{\THH{Characterizing functions ($\cF_i\in \mathfrak{F}$)}}
In order to characterize the system at a microscopic level, one may choose to equip $\MM$ with a finite set of characterizing functions, labeled by $i\in\II$, denoted $\cF_{i}:\MM\rightarrow\FF_i$.

\begin{definition}[Extensivity]
A function $\cF_i$ is extensive if $\FF_i$ is a $\RR^+$-semimodule (see Def.~\ref{def:module}) and if it verifies
\begin{equation}
\bm{m}_3 =\bm{m}_1\boldsymbol{\cdot} \bm{m}_2 \Rightarrow \cF_i(\bm{m}_3) = \cF_i(\bm{m}_1) + \cF_i(\bm{m}_2) \quad \forall (\bm{m}_1, \bm{m}_2, \bm{m}_3) \in \MM^3.
\label{eq:extensivity}
\end{equation}
\end{definition}

\begin{definition}[$\RR^+$-semimodule]
A set $\mathbb{S}$ is a $\RR^+$-semimodule if for all $(r_1, r_2) \in {\RR^+}^2$ and $(s_1, s_2) \in \mathbb{S}^2$, the following properties hold
\begin{enumerate}
\item $r_1\,(s_1 + s_2) = r_1\,s_1 + r_2\,s_2$,
\item $(r_1 + r_2)\,s_1 = r_1\,s_1 + r_2\,s_1$,
\item $(r_1\,r_2)\,s_1 = r_1\,(r_2\,s_1)$,
\item $1\,s_1 = s_1$,
\item $0\, s_1 = 0$.
\end{enumerate}
\label{def:module}
\end{definition}

\setcounter{example}{0}
\begin{example}[Continued]
For a gas of $N$ identical, non-interacting particles, the function $\cF_e : \MM \mapsto \RR^+$ defined as 
\begin{equation}
\cF_e : \bm{m} \longmapsto \cF_e(\bm{m}) = \sum_{i = 1}^{N} \frac{\norm{\bm{p}^i(\bm{m})}^2}{2\,m}
\end{equation}
gives the energy of the system in microstate $\bm{m}$, where $\bm{p}^i(\bm{m})$ is the momentum of particle $i$, and $m$ is the mass of a particle. It fulfills the extensivity property defined in Eq.~(\ref{eq:extensivity}).
\end{example}
\setcounter{example}{1}
\begin{example}[Continued]
In the Ising model~\THH{\cite{liechtenstein1984exchange}}, the function $\cF_e : \MM \mapsto \RR$ defined as 
\begin{equation}
\cF_e : \bm{m} \longmapsto \cF_e(\bm{m}) = -\frac{1}{2}\,\bm{m}^\intercal \,\bm{J}_\mathrm{ex}\,\bm{m}
\end{equation}
gives the energy of the system in microstate $\bm{m}$, where each coefficient ${\bm{J}_\mathrm{ex}}_{i,j}$ is the exchange energy between atom $i$ and atom $j$. It does not fulfill the extensivity property defined in Eq.~(\ref{eq:extensivity}).
\end{example}
\begin{example}
The function $\cF_n : \MM \mapsto \mathbb{N}^+$ defined as
\begin{equation}
\cF_n : \bm{m} \longmapsto \cF_n(\bm{m})
\end{equation}
where $\cF_n(\bm{m})$ is the number of particles of the system in microstate $\bm{m}$ \THH{is extensive}.
\end{example}

\begin{example}\label{ex:cFr}
The function $\cF_r$ defined as
\begin{equation}
\cF_r(\bm{m}) = \left(\bm{r}^i\right)_{1\leq i \leq \cF_n(\bm{m})}, 
\end{equation}
where $\bm{r}^i \in \RR^3$ is the position of particle $i$, gives the set of all particle positions for the system in microstate $\bm{m}$. It is not extensive.
\end{example}

\begin{example}\label{ex:Fv}
The function $\cF_v : \MM \mapsto \RR^+$ defined as 
\begin{equation}
\cF_v : \bm{m} \longmapsto \cF_v(\bm{m})
\end{equation}
where $\cF_v(\bm{m})$ is the volume occupied by the system in microstate $\bm{m}$.
\THH{The choice} of such function $\cF_v$ is hardly unique 
\THH{(see remark~\ref{rk:Volume1} below).}
Here, we propose to define $\cF_v(\bm{m})$ as the minimal bounding volume enclosing all particle positions of microstate $\bm{m}$ that accounts for the container geometry and its degrees of freedom.
For instance, for a cylindrical container of fixed base $A$ closed by a piston moving freely along axis $z$, we can define the volume as
\begin{equation}
\cF_v(\bm{m}) = A \times h(\bm{m}),\quad \text{with } h(\bm{m}) = \max\{r^i_z\,|\, \bm{r}^i \in \cF_r(\bm{m})\}.
\end{equation}
\THH{This function} does not fulfill the extensivity property defined in Eq.~(\ref{eq:extensivity}).
\end{example}

\THH{
\begin{remark}[Volume]\label{rk:Volume1}
The mathematical conceptualization of the volume is  a challenging issue\footnote{In~\cite[Chap.\,3,p.3]{Grothendieck1986}, Grothendieck mentions the absence (in most textbooks) of any "serious" definition of the notion of length (of a curve), of area (of a surface), of volume (of a solid).}.
In the context of SP, its physical conceptualization is also an issue:
\begin{itemize}
\item[(i)] a possible choice could be the volume occupied by the particles, e.g. that delimited by the 3D simplicial envelope of all the particle positions,
\item[(ii)] an alternative is to consider the volume of a container, in which the particle are authorized to evolve.
\end{itemize}
The case~(i) allows the definition of a characterizing function $\cF_v$, the volume being intrinsically related to the microstate.
In (ii), the microstate does not encode the information of the container volume: its set of particles in contact with the container boundary can even be empty. This information is an "experimental constraint" (presented in section~\ref{sec:Ma} below).
Note also that example~\ref{ex:Fv} corresponds to a hybrid description in between (i) and (ii).
\end{remark}
}

In the following, we denote $\mathfrak{F} = \{\cF_i\}_{i \in \II}$ the set of characterizing functions on $\MM$.

\section{Experimental conditions and accessible microstates \THH{
($\textbf{m}\in\MM_a$)
}}
\label{sec:Ma}
Experimental conditions may constrain characterizing functions to take values that are compatible with these experimental conditions. Thus, under experimental conditions, the configuration space becomes restricted to a set of \emph{accessible} microstates $\MM_a  \subset \MM$.
\begin{definition}[Set of accessible microstates $\MM_a$]
Denote $\II^\mathrm{fixed}\subseteq\II$ the set of labels of characterizing functions that are experimentally constrained. Due to the constraints, a function $\cF_i$, $i \in \II^\mathrm{fixed}$ can only take admissible values in $\FF_i^\mathrm{fixed} \subseteq \FF_i$.

Denote $\bm{\theta}^\mathrm{fixed} :=\left(\FF_i^\mathrm{fixed}\right)_{i\in\II^\mathrm{fixed}}$. The set of accessible microstates $\MM_a\left(\bm{\theta}^\mathrm{fixed}\right)$ is
    \begin{equation}
        \MMa\big( \bm{\theta}^\mathrm{fixed} \big) = \left\{ \bm m \in \MM \;|\; \cF_i(\bm{m}) \in \FF_i^\mathrm{fixed} \quad \forall i \in  \II^\mathrm{fixed}\right\}.
    \end{equation}
 \end{definition}
Remark that $\II^\mathrm{fixed}=\II$ defines an isolated system with respect to the chosen characterizing functions.

\begin{figure}[tb]
\captionsetup{width= \textwidth}
\centering
\begin{subfigure}[t]{0.29 \textwidth}
\centering
\begin{tikzpicture}[scale = 0.8]
\node [circle, draw, red] (a) at (0, 0) {};
\node [circle, draw, red] (b) at (1, 1) {};
\node [circle, draw, red] (c) at (0.2, 2) {};
\draw (-0.5, -0.5) [thick, rounded corners = 5pt] rectangle (1.5, 2.5);
\end{tikzpicture}
\caption{Fixed number of particles and fixed volume.}
\end{subfigure}
\begin{subfigure}[t]{0.29 \textwidth}
\centering
\begin{tikzpicture}[scale = 0.8]
\node [circle, draw, red] (a) at (0, 0) {};
\node [circle, draw, red] (b) at (1, 0) {};
\node [circle, draw, red] (c) at (0.3, 0.6) {};
\draw (-0.5, -0.5) [thick, rounded corners = 5pt] rectangle (1.5, 2.5);
\draw (-0.5, 2.3) [very thick, fill = white, white, opacity = 1] rectangle (1.5, 3);
\draw (-0.3, 1) [thick] -- ++(1.6, 0);
\draw (0.5, 1) [thick] -- ++(0, 2);
\end{tikzpicture}
\caption{Fixed number of particles.}
\end{subfigure}
\begin{subfigure}[t]{0.29\textwidth}
\centering
\begin{tikzpicture}[scale = 0.8]
\node [circle, draw, red] (a) at (0, 0) {};
\node [circle, draw, red] (b) at (1, 1) {};
\draw (-0.5, -0.5) [thick, rounded corners = 5pt] rectangle (1.5, 2.5);
\draw (-0.5, 2.3) [very thick, fill = white, white, opacity = 1] rectangle (1.5, 3);
\draw (-0.5, 2.3) [thick, dashed] -- ++ (2, 0);
\node [circle, draw, red] (c) at (1, 3) {};
\end{tikzpicture}
\caption{Fixed volume.}
\end{subfigure}
\caption{Examples of experimental conditions for a gas in a tank.}
\label{fig:thermo_constraints}
\end{figure}
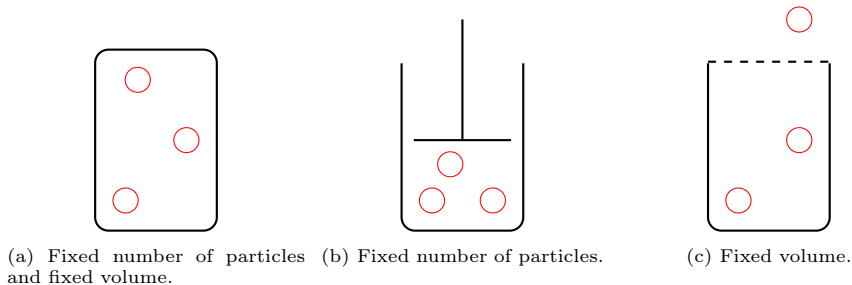

\setcounter{example}{0}
\begin{example}[Continued]
Consider a gas of $N$ particles in a closed tank. The system cannot exchange particles with the environment, therefore the number of particles $\cF_n(\bm{m})$ is fixed to $N$. Denote $\FF_n^\mathrm{fixed} = \{N\}$. The set of accessible microstates is $\MM_a\left(N\right) = \left\{\bm{m} \in \mathbb{M} \;|\; \cF_n(\bm{m}) \in \FF_n^\mathrm{fixed}\right\}$.
\end{example} 
\setcounter{example}{0}
\begin{example}[Continued]
Consider a gas of $N$ particles in a closed tank occupying a space $\Pi \subset \RR^3$. Denoting $\FF_n^\mathrm{fixed} = \{N\}$ and $\FF_r^\mathrm{fixed} = \Pi^N$, the set of accessible microstates is
$
\MM_a\left(N,\, \Pi\right) = \left\{\bm{m} \in \MM \;|\; \cF_n(\bm{m})\in \FF_n^\mathrm{fixed},\, \cF_r(\bm{m}) \in \FF_r^\mathrm{fixed}\right\}
$.
\THH{Note that the constraint on $\FF_r^\mathrm{fixed}$ corresponds to the case (ii) in remark~\ref{rk:Volume1}, the container being described by $\Pi$.}
\end{example}
Other examples of experimental conditions are shown on Fig.~\ref{fig:thermo_constraints}.

\THH{
\paragraph{Intermediate point}
At this step, the physics is described through:
$\MM$ (microscopic representation of physical particles, i.e. microstates),
$\MM_a\subseteq\MM$ (microstates accessible under fixed macroscopic experimental constraints),
$\mathfrak{F}$ (links to quantities that can be observed at the macroscopic scale).
This is not sufficient to derive an autonomous physical law at macroscopic scale.
This issue can be solved by: (i)~completing $\mathfrak{F}$ with a single new function, namely, the surprisal, which is elaborated from a stochastic description, giving rise to the entropy; and (ii)~applying the fundamental principle introduced by Boltzmann~\cite{boltzmann1877beziehung} to derive a microstate probability that physically makes sense.

The following sections revisit this approach, focusing on the propensity of macroscopic quantities to communicate with their peer in an external environment.
To this end,} we assume that all characterizing functions that are not explicitly fixed by experimental conditions can still depend on microstate $\bm{m}$, and we denote $\II^\mathrm{free} := \II \backslash \II^\mathrm{fixed}$ the set of labels of characterizing functions not fixed by experimental conditions.


\section{Stochastic representation and measure of uncertainty}
\label{sec:stoch}

\subsection{Microstate stochastic description}
The system fluctuates from one accessible microstate to another. 
\THH{It is considered to be impossible to predict these fluctuations in a deterministic fashion at the macroscopic level:  \ac{sp} adopts a stochastic framework that model their random description.
}
Indeed, from Prop.~(\ref{prop:measure}), $\left(\MM, \mathcal{P}(\MM)\right)$ is measurable, therefore so is $\left(\MM_a, \text{Tr}\left(\mathcal{P}(\MM)\right)_{\MMa}\right)$, where $\text{Tr}\left(\mathcal{P}(\MM)\right)_{\MMa}$ denotes the trace of $\mathcal{P}(\MM)$ on $\MM_a$~\cite{plachky2001ideal}. 
\THH{Assuming that $\MM_a$ is countable and that the distribution $p$ is discrete\footnote{\THH{
Extensions to continuous measurable spaces are available and similar, replacing the sum by an integral with a Lebesgues measure and using distributions in (\ref{eq:exp}).
}},
}
we can define a probability distribution 
\THH{
$p\left(.\, | \,\bm{\theta}^\mathrm{fixed}\right) :\MM_a\left(\bm{\theta}^\mathrm{fixed}\right) \mapsto [0, 1]$, denoted $p$ for short below,} which
assigns to each microstate $\bm{m} \in \MM_a\left(\bm{\theta}^\mathrm{fixed}\right)$ a probability $p(\bm{m})$ to be the actual microstate of the system \THH{(the Boltzmann principle in section~\ref{sec:proba} will provide a tool to determine this probability).}

\THH{The} average of a random quantity $\cF(\bm{m})$ is given by its expectation $\EE_p[\cF]$, defined as
\begin{equation}
\EE_p[\cF] = \sum_{\bm{m}\in \MM_a\left(\bm{\theta}^\mathrm{fixed}\right)}p(\bm{m})\,\cF(\bm{m}).
\label{eq:exp}
\end{equation}

\subsection{Statistical entropy}
Given some basis of \THH{information units}  $b > 1$, a microstate $\bm{m}$ with probability $p(\bm{m})$
has a surprisal $\cS_p^b(\bm{m})$ defined as
\begin{equation}
\cS_p^b(\bm{m}) = \log_b \frac{1}{p(\bm{m})} \quad \text{with } \THH{\log_b = \frac{1}{\ln(b)}\,\ln}.
\end{equation}
The surprisal, or information content, quantifies how much the occurrence of microstate $\bm{m}$ is surprising. For example, if some microstate $\bm{m}$ is the state of the system for certain, it has probability 1 and surprisal 0.

As the probability of two independent events $\bm{m}_1$ and $\bm{m}_2$ verifies 
\begin{equation}
p\left(\bm{m}_1 \boldsymbol{\cdot}\bm{m}_2\right) = p(\bm{m}_1)\,p(\bm{m}_2),
\end{equation}
the surprisal function $\cS_p^b$ verifies the extensivity property defined in Eq.~(\ref{eq:extensivity}).

The surprisal allows the definition of a measure of lack of information on average for a probability distribution $p$ and a basis $b$, namely, the statistical entropy $\mathsf{S}_b(p)$~\cite{gray2011entropy} defined as
\begin{equation}
\mathsf{S}_b(p) =  \EE_p[\cS_p^b].
\label{eq:entropy}
\end{equation}
The statistical entropy can be interpreted of as ``the average number of questions to ask with $b$ possible answers per question'' in order to know the actual microstate for certain.

\begin{figure}[tb]
\centering
\begin{subfigure}[b]{0.2\textwidth}
\centering
\begin{tikzpicture}
\draw (0, 0) node[circle, draw, fill = red, red]{};
\draw (0.5, 0) node[circle, draw, fill = red, red]{};
\end{tikzpicture}
\caption{}
\end{subfigure}
\begin{subfigure}[b]{0.2\textwidth}
\centering
\begin{tikzpicture}
\draw (0, 0) node[circle, draw, fill = red, red]{};
\draw (0.5, 0) node[circle, draw, fill = blue, blue]{};
\end{tikzpicture}
\caption{}
\end{subfigure}
\begin{subfigure}[b]{0.2\textwidth}
\centering
\begin{tikzpicture}
\draw (0, 0) node[circle, draw, fill = blue, blue]{};
\draw (0.5, 0) node[circle, draw, fill = blue, blue]{};
\end{tikzpicture}
\caption{}
\end{subfigure}
\begin{subfigure}[b]{0.2\textwidth}
\centering
\begin{tikzpicture}
\draw (0, 0) node[circle, draw, fill = blue, blue]{};
\draw (0.5, 0) node[circle, draw, fill = red, red]{};
\end{tikzpicture}
\caption{}
\end{subfigure}
\caption{Possible outcomes for a coin tossed twice.}
\label{fig:statistical_entropy}
\end{figure}
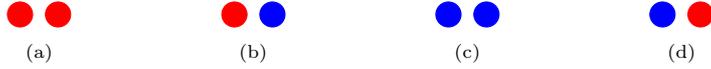

\setcounter{example}{5}
\begin{example}
Consider the outcomes of tossing a coin twice. The coin can come up heads or tails after each toss, hence $2 \times 2 = 4$ possible outcomes (Fig.~\ref{fig:statistical_entropy}).
If all outcomes are equiprobable, one needs at least two questions with two possible answers each to know the exact outcome: 
\begin{enumerate}
\item Did the coin come up heads or tails after the first toss? 
\item Did the coin come up heads or tails after the second toss?
\end{enumerate} 
As it happens, taking $p_1:\, \bm{m} \mapsto p_1(\bm{m}) = \frac{1}{4}$ and $b$ = 2 in Eq.~(\ref{eq:entropy}) yields $\mathsf{S}_b(p_1) = -\log_2\frac{1}{4} = 2$.

However, if the probability distribution is not uniform, some outcomes are more probable than others, and the uncertainty is lower; ditto the entropy. For instance, with a probability distribution $p_2$ assigning $\frac{1}{2}$ to outcome (A), $\frac{1}{4}$ to outcome (B), and $\frac{1}{8}$ to outcomes (C) and (D), the entropy becomes $\mathsf{S}_b(p_2) = 1.75 < \mathsf{S}_b(p_1) = 2$.
\label{ex:coin}
\end{example}
In information theory, statistical entropy relates to optimal encoding of information. Suppose you repeat the coin toss experiment of Ex.~(\ref{ex:coin}) for a long period of time, and wish to record every outcome on a computer. For a sequence of two tosses with distribution $p_1$, an outcome cannot be encoded in less than two bits; while with distribution $p_2$, outcome (A) can be encoded on one bit, outcome (B) on two, and outcomes (C) and (D) on three, that is, $1 \times \frac{1}{2} + 2 \times \frac{1}{4} + 6 \times \frac{1}{8} = 1.75$ bits on average. The most frequent outcome takes the least encoding space; conversely, the comparatively large encoding space taken by outcomes (C) and (D) is compensated by the rarity of their occurrence. On the whole, exploiting the knowledge underpinned by distribution $p_2$ reduces the encoding cost. This principle underlies Morse code (and, more generally, lossless entropy encoding like Huffman coding~\cite{knuth1985dynamic}): very common letters such as ``e'' or ``i'' take much fewer dots than less common letters like ``j'' or ``q''.

For compacity, the statistical entropy becomes in the following
\begin{equation}
\mathsf{S}^k(p) := \mathsf{S}_{b = \exp(\nicefrac{1}{k})}(p) = -k\,\sum_{\bm{m} \in \MM_a(\bm{\theta}^\mathrm{fixed})} p(\bm{m})\ln p(\bm{m}),
\label{eq:entropy_k}
\end{equation}
for all $p$ defined on $\MM_a(\bm{\theta}^\mathrm{fixed})$.

\THH{
\begin{remark}[Random structures]\label{rk:combinatorics2}
Following remark~\ref{rk:combinatorics1}, the cases of microstates involving combinatorial structures based on elaborated specifications (such as molecules and chemical reaction processes) require elaborate tools addressing random structures (see e.g.~\cite[part.\,C]{flajolet2009analytic}) that are out of the scope of this paper.
\end{remark}
}


\section{Microstate probability distribution at equilibrium and partition function}
\label{sec:proba}
\subsection{Thermodynamic equilibrium}
A system is at thermodynamic equilibrium when its statistics stops evolving. At this point, the \emph{ergodic hypothesis} postulates that over a ``sufficiently long'' period of time $t$, the system explores all its accessible microstates.
Assuming that a microstate $\bm{m}$ can be measured at a time $\tau$ through
$\mathcal{M} :\, \RR^+ \mapsto \MM_a$, this means that at thermodynamic equilibrium, the temporal mean of a quantity $\overline{\cF}_i$ coincides with its expectation $\EE_p [\cF_i]$: 
\begin{equation}
\overline{\cF}_i := \lim_{t\rightarrow +\infty}\frac{1}{t}\int_0^t\del{\cF_i\circ\mathcal{M}}(\tau)\,\text{d}\tau = \EE_p [\cF_i]. 
\label{eq:ergodicity}
\end{equation}

Therefore, for a given set $\bm{\theta}^\mathrm{free}$ of mean values $(\overline{\cF}_i)_{i \in \II^\mathrm{free}}$, the ergodic hypothesis translates into a set of hypotheses $\mathsf{H}\left(\bm{\theta}^\mathrm{free}\right)$ defined as
\begin{equation}
\mathsf{H}\left(\bm{\theta}^\mathrm{free}\right) =\left(\EE_p[\cF_i] = \overline{\cF}_i\right)_{i \in \II^\mathrm{free}}.
\label{eq:H}
\end{equation}

While still discussed~\cite{patrascioiu1987ergodic} (especially regarding the definition of ``sufficiently long''), this hypothesis is the foundation of equilibrium statistical physics, and we assume its validity in the following.

\subsection{\THH{Boltzmann principle: maximum entropy at thermodynamic equilibrium (reminder)}}

By definition, a system at equilibrium does not evolve. Since change is new information, the information given by a system at equilibrium is minimal.
As statistical entropy is a measure of lack of information, 
it follows that at equilibrium, the entropy is maximal: this is the Boltzmann principle.

It follows that the microstate probability distribution at thermodynamic equilibrium $p^\star$ is
\begin{align}
p^\star 
& = \underset{p}{\arg\max} \;\mathsf{S}^k (p)
\text{~~subject to} \hspace*{-8pt}
\sum_{\bm{m} \in \MM_a\left(\bm{\theta}^\mathrm{fixed}\right)}p(\bm{m}) = 1
\nonumber\\
& \THH{= \underset{p}{\arg}\;\;\underset{p,\lambda}{\max} \;\mathsf{S}^k (p) + \lambda\,\left( \sum_{\bm{m} \in \MM_a\left(\bm{\theta}^\mathrm{fixed}\right)}p(\bm{m}) - 1 \right),}
\label{eq:entropy_max1}
\end{align}
where the second line specifies the constraint using a Lagrange multiplier $\lambda$.

\subsection{\THH{
Boltzmann principle in a macroscopic external environment}}
\THH{
The macroscopic quantities ($\cF_{i\in \II^\mathrm{free}}(\bm{m})$) whose fluctuation is experimentally allowed are prone to communicate with the external environment (for example, through exchanges of particle, energy, etc).
Assuming ergodicity and an infinite ratio between macroscopic and microscopic time scales, this means that the expectation of these fluctuating quantities (with value $\overline{\cF}_{i\in \II^\mathrm{free}}$) can change over time (at the slow macroscopic scale), while the thermodynamic equilibrium is satisfied (at the fast microscopic scale) and continuously updated (at the slow macroscopic scale).

As a main step of this paper, this issue is addressed by deriving the conditional probability $p^\star\left(\bm{m}\;|\;\mathsf{H}\left(\bm{\theta}^\mathrm{free}\right)\right)$ that 
maximizes entropy (Boltzmann principle) constrained by
the given macroscopic values $\bm{\theta}^\mathrm{free}=(\overline{\cF}_i)_{i \in \II^\mathrm{free}}$, namely,
}
\begin{align}
p^\star = &\underset{p}{\arg\max} \;\mathsf{S}^k (p) 
\text{\small~~subject to~}
\begin{dcases}
\sum_{\bm{m} \in \MM_a\left(\bm{\theta}^\mathrm{fixed}\right)}p(\bm{m}) = 1,\\
\mathsf{H}\left(\bm{\theta}^\mathrm{free}\right),
\end{dcases}
\nonumber\\
= & \THH{
\underset{p}{\arg}\;\;\underset{p,\lambda_0,\lambda_{i\in \II^\mathrm{free}}\!\!\!}{\max} \;\mathsf{S}^k (p) + \lambda_0\,\left( \sum_{\bm{m} \in \MM_a\left(\bm{\theta}^\mathrm{fixed}\right)}p(\bm{m}) - 1 \right)
}
\nonumber\\
& \THH{+ \sum_{i \in \II^\mathrm{free}}
\lambda_i \, \left( \EE_p[\cF_i] - \overline{\cF}_i\right),}
\label{eq:entropy_max}
\end{align}
where $\mathsf{H}\left(\bm{\theta}^\mathrm{free}\right)$ accounts for the experimental conditions (see Eq.~(\ref{eq:H})).

Note that, by definition, this probability naturally restores that of (\ref{eq:entropy_max1}), if 
the values $\overline{\cF}_{i\in \II^\mathrm{free}}$ are the expectations of ${\cF}_{i\in \II^\mathrm{free}}$ computed for 
probability (\ref{eq:entropy_max1}).

\begin{theorem}
\label{th:boltz}
Let $\bm{\theta}^\mathrm{free} := \left(\overline{\cF}_i\right)_{i\in \II^\mathrm{free}} \in  \underset{i\in \II^\mathrm{free}}{\text{\Large $\times$}} \FF_i$, where $\underset{i\in \II^\mathrm{free}}{\text{\Large $\times$}} \FF_i$ denotes the Cartesian product of the $\left(\FF_i\right)_{i \in \II^\mathrm{free}}$.

Then for all $\bm{m} \in \MM_a \left(\bm{\theta}^\mathrm{fixed}\right)$,
\begin{subequations}
\begin{equation}
p^\star\left(\bm{m}\;|\;\mathsf{H}\left(\bm{\theta}^\mathrm{free}\right)\right) = \frac{\exp\del{\frac{\sum_{i\in \II^\mathrm{free}}\lambda_i\,\cF_i(\bm{m})}{k}}} {\cZ\left( \bm{\lambda}^\mathrm{free}\right)},
\label{eq:p_opt}
\end{equation}
where
\begin{equation}
\cZ\left(\bm{\lambda}^\mathrm{free}\right) := \sum_{\bm{m}\in\MM_a\left(\bm{\theta}^\mathrm{fixed}\right)}\exp\del{\frac{\sum_{i\in \II^\mathrm{free}}\lambda_i\,\cF_i(\bm{m})}{k}}
\end{equation}
is the partition function of the system,
and, for all $i \in \II^\mathrm{free}$, $\lambda_i$ verifies
\begin{equation}
\dpd{}{\lambda_i} k\,\ln\cZ\left(\bm{\lambda}^\mathrm{free}\right) = \overline{\cF}_i.
\end{equation}
\end{subequations}
\end{theorem}
The proof is in \ref{app:proof_boltzmann}.

\begin{definition}[Thermodynamic entropy]
The thermodynamic entropy $S^k\left(\bm{\theta}^\mathrm{free}\right)$ is defined as the statistical entropy for the probability distribution at equilibrium given $\bm{\theta}^\mathrm{free}$:
\begin{equation}
S^k\left(\bm{\theta}^\mathrm{free}\right) = \mathsf{S}^k\left(p^\star\left(\cdot \;|\; \mathsf{H}\left(\bm{\theta}^\mathrm{free}\right)\right)\right).
\label{eq:entropy_thermo}
\end{equation}
\end{definition}
\begin{property}
The thermodynamic entropy function $S^k$ is a Legendre transform of $k\,\ln\cZ$ and we have
\begin{equation}
S^k\left(\bm{\theta}^\mathrm{free}\right) = k\,\ln\cZ\left( \bm{\lambda}^\mathrm{free}\right) -\sum_{i\in \II^\mathrm{free}}\lambda_i\,\overline{\cF}_i.
\label{eq:entropy_legendre}
\end{equation}
\begin{proof}
\begin{equation*}
\begin{split}
S^k\left(\bm{\theta}^\mathrm{free}\right) &\overset{(a)}{=} \mathsf{S}^k\left(p^\star\left(\cdot \;|\; \mathsf{H}\left(\bm{\theta}^\mathrm{free}\right)\right)\right)\\
&\overset{(b)}{=} -k\sum_{\bm{m}\in\MM_a\left(\bm{\theta}^\mathrm{fixed}\right)}p^\star\left(\bm{m}\;|\; \mathsf{H}\left(\bm{\theta}^\mathrm{free}\right)\right)\,\ln p^\star\left(\bm{m}\;|\; \mathsf{H}\left(\bm{\theta}^\mathrm{free}\right)\right)\\
&\overset{(c)}{=} -k\sum_{\bm{m}\in\MM_a\left(\bm{\theta}^\mathrm{fixed}\right)}p^\star\left(\bm{m}\;|\; \mathsf{H}\left(\bm{\theta}^\mathrm{free}\right)\right)\,
\ln\del{\frac{\exp\del{\frac{\sum_{i\in \II^\mathrm{free}}\lambda_i\,\cF_i(\bm{m})}{k}}} {\cZ\left( \bm{\lambda}^\mathrm{free}\right)}}\\
&= -k\sum_{\bm{m}\in\MM_a\left(\bm{\theta}^\mathrm{fixed}\right)}p^\star\left(\bm{m}\;|\; \mathsf{H}\left(\bm{\theta}^\mathrm{free}\right)\right)\del{\frac{\sum_{i\in \II^\mathrm{free}}\lambda_i\,\cF_i(\bm{m})}{k} - \ln\cZ\left( \bm{\lambda}^\mathrm{free}\right)}\\
&\overset{(d)}{=} -\sum_{i\in \II^\mathrm{free}}\lambda_i\,\overline{\cF}_i + k\,\ln\cZ\left( \bm{\lambda}^\mathrm{free}\right),
\end{split}
\end{equation*}
using (a) Eq.~(\ref{eq:entropy_thermo}), (b) Eq.~(\ref{eq:entropy_k}), (c) Eq.~(\ref{eq:p_opt}), and (d) Eqs.~(\ref{eq:exp})-(\ref{eq:ergodicity}).

We deduce that $S^k$ is a Legendre transform  
of $k\,\ln\cZ$ (see also \cite{zia2009making}).
\end{proof}
\label{prop:legendre}
\end{property}
\begin{property}
It follows from Prop.~(\ref{prop:legendre}) that for all $i \in \II^\mathrm{free}$, the Lagrange multiplier $\lambda_i$ is the derivative of the thermodynamic entropy function with respect to average $\overline{\cF}_i$
\begin{equation}
\lambda_i = - \dpd{S^k}{\overline{\cF}_i}\left(\bm{\theta}^\mathrm{free}\right).
\label{eq:lambda}
\end{equation}
\end{property}
\begin{example}
In particular, this defines the system temperature $T$, chemical potential $\mu$, and pressure $P$ as
\begin{align}
\frac{1}{T} &:= \dpd{S^k}{\overline{\cF}_e}\left(\bm{\theta}^\mathrm{free}\right) , &
\frac{\mu}{T} &:= -\dpd{S^k}{\overline{\cF}_n}\left(\bm{\theta}^\mathrm{free}\right), & \frac{P}{T} &:=\dpd{S^k}{\overline{\cF}_v}\left(\bm{\theta}^\mathrm{free}\right).
\label{eq:T}
\end{align}
\end{example}
\paragraph{Adiabatic process}
For a system going through an adiabatic process (no thermal exchange with the environment), the surprisal is independent of $\bm{m}$ so that
\begin{equation}
\cS_{p^\star}^b(\bm{m}) = S \quad \forall \bm{m}.
\end{equation}
That implies that for such systems, all microstates have the same probability
\begin{equation}
p^\star(\bm{m}) = \frac{1}{\Omega}, \quad \text{with }\Omega = \text{card}\del{\MM_a\left(\bm{\theta}^\mathrm{fixed}\right)}.
\label{eq:equiprob}
\end{equation}
From Eq.~(\ref{eq:p_opt}), it follows that for such a system, we have
\begin{equation}
\sum_{i\in \II^\mathrm{free}}\lambda_i\,\cF_i(\bm{m}) = C,
\end{equation}
where $C$ is independent of $\bm{m}$.
\begin{example}
In particular, a system that can only exchange volume with its environment verifies $\cF_e(\bm{m}) + P\,\cF_v(\bm{m}) = C$, where $C$ is independent of $\bm{m}$.
\end{example}

\subsection{Identification of Boltzmann constant}
To ensure that the statistical entropy does coincide with the thermodynamic entropy at equilibrium, the constant $k$ must be chosen as the Boltzmann constant\footnote{\THH{
Note that, from definition (\ref{eq:entropy_k}), choosing the unit USI reference $k_0=1\,$J.K$^{-1}$ as a unit information quantity,
this value corresponds to the question number base $b_B=\exp(k_0/k_B)= \exp(10^{23}/1.38) \approx 10^{3.147\mathrm{e}+22}\approx 2^{9.4736\mathrm{e}+21}$.
This means that explaining +1\,J.K$^{-1}$ requires  about $10^{22}$~bits for a gas.
}}
$k_B = 1.38 \times 10^{-23}$ J.K$^{-1}$.
Indeed, consider an ideal gas of $N$ non-interacting atoms in a box of volume $V$ at temperature $T$, represented by their position and momentum. The partition function $\cZ$ is given by
\begin{equation}
\cZ(T \,|\,N, V) = V^N\del{\frac{2\,\pi\,m \,k \, T}{h ^2}}^{\nicefrac{3\,N}{2}},
\end{equation}
where here $m$ denotes the mass of an atom, and $h$ is the Planck constant.
From Prop.~(\ref{prop:legendre}), the thermodynamic entropy $S^k(\overline \cF_e, N, V)$ is given by
\begin{equation}
S^k(\overline \cF_e, N, V) = k\,\ln \cZ(T\,|\,N, V) + \frac{\overline \cF_e}{T}. 
\end{equation}
Moreover, from Eq.~(\ref{eq:T}), the pressure $P$ is given by
\begin{equation}
P = T\,\dpd{S^k}{V}(\overline \cF_e, N, V) = \frac{N\,k\,T}{V}.
\end{equation}
Therefore, $k$ must be identified with $k_B$ so that the ideal gas law $P\,V = N\,k_B\,T$ is verified.

\section{Final \ac{phs} model}\label{sec:final}
\THH{
The macrosopic description of open system can be achieved by using balanced equations of variations of entropy, energy, mass, etc.
Port-Hamiltonian systems provide an adapted framework for such physical descriptions.
This section addresses this issue by using a standard formulation~(recalled in section~\ref{ssec:PHSreminder}), in which the energy is expressed as a function of entropy and state variables: this requires to invert $S:(E,\dots)\mapsto S(E,\dots)$ w.r.t. $E$, to introduce the hamiltonian $H=E:(S,\dots)\mapsto H(S,\dots)$ in a first step.
Note that this inversion is a technical step that could be avoided (to still use the entropy function) by considering contact forms as in~\cite{van2018geometry} (see also \cite{ottinger2005beyond} for the alternative GENERIC formulation).
}

\subsection{Reminder on port-Hamiltonian systems}\label{ssec:PHSreminder}
The \ac{phs} formalism provides a unified formalism for the modeling of multiphysical systems, in the sense that it recognizes energy as a universal currency. Indeed, any physical system can be divided into parts that interact with each other via energy exchanges.

Detailed presentations of \ac{phs} are available in~\cite{duindam2009modeling, van2014port}. In this paper, we rely on a differential-algebraic formulation adapted to multiphysical systems~\cite{falaize2016passive, remy2021time}.
This formulation allows the representation of a dynamical system as a network of
\begin{enumerate}
\item storage components of state $\bm{x}$ and energy $E(\bm{x})$;
\item passive memoryless components described by an effort law $z : \bm{w}\mapsto z(\bm{w})$, such as the dissipated power $P_\mathrm{d}=z\del{\bm{w}}^\intercal\bm{w}$ is non-negative for all flows $\bm{w}$;
\item connection ports conveying the \emph{outgoing} power 
$P_\mathrm{ext} = \bm{u}^\intercal \bm{y}$ where $\bm{u}$ are inputs and $\bm{y}$ are outputs.
\end{enumerate}
The system flows $\bm{f}$ and efforts $\bm{e}$ are coupled through a (possibly dependent on $\bm{x}$) skew-symmetric interconnection matrix $\bm{S} = -\bm{S}^\intercal$, so that
\begin{equation}
				\underbrace{\begin{bmatrix}
					\dot {\bm{x}}\\
					\bm{w}\\
					\bm{y}
				\end{bmatrix}}_{\bm{f}}
 				= \bm{S}
 				\underbrace{\begin{bmatrix}
					\nabla E(\bm{x})\\
					z(\bm{w})\\
					\bm{u}
 				\end{bmatrix}}_{\bm{e}}.
\label{eq:SHP}
\end{equation}

Such systems satisfy the power balance 
\begin{equation}
P_\mathrm{s} + P_\mathrm{d} + P_\mathrm{ext} = 0
\label{eq:powerbalance}
\end{equation}
where $P_\mathrm{s} = \nabla E(\bm{x})^\intercal\,\dot{\bm{x}}$ denotes the stored power.
\begin{proof}
\begin{equation*}
P_\mathrm{s} + P_\mathrm{d} + P_\mathrm{ext} = \nabla E(\bm{x})^\intercal\,\dot{\bm{x}} + z\del{\bm{w}}^\intercal\,\bm{w} + \bm{u}^\intercal\,\bm{y} = \bm{e}^\intercal \bm{f} = \bm{e}^\intercal \bm{S} \bm{e} =(\bm{e}^\intercal \bm{S}\bm{e})^\intercal= -\bm{e}^\intercal \bm{S}\bm{e} = 0
\end{equation*} 
due to the skew-symmetry of $\bm{S}$. 
\end{proof}

Note that in this paper, we adopt the \emph{passive sign convention} (also called receiver convention) for all components, including external sources. This means that a flow is defined positive when entering the component ~\cite{bigelow2020power}.

\subsection{Macroscopic state and energy}
In the previous sections, we derived the thermodynamic entropy as a function of the energy and other macroscopic variables. In order to obtain a port-Hamiltonian formulation such as Eq.~(\ref{eq:SHP}), we choose to express the energy as a function of the thermodynamic entropy and other macroscopic variables instead.

Denote $\overline \cS := S^k\left(\bm{\theta}^\mathrm{free}\right)$, and $\bm{\theta}^x := \bm{\theta}^\mathrm{free}\backslash \overline \cF_e$ the set of macroscopic quantities that are not the energy, with corresponding set of labels $\II^x$.
We choose to define the (extensive) macroscopic state $\bm{x}$ as
\begin{equation}
\bm{x} = \left[\overline{\cS}, \bm{\theta}^x\right]^\intercal,
\end{equation}
so that the flow $\dot{\bm{x}}$ accounts
for the time variation of extensive quantities.
Assuming that the entropy function $S^k$ is invertible with respect to $\overline{\cF_e}$, we
define the macroscopic energy function $E$ as
\begin{equation}
E : \bm{x} \longmapsto E\left(S^k(\overline \cF_e,\, \bm{\theta}^x),\, \bm{\theta}^x\right) = \overline \cF_e,
\end{equation}
so that the effort $\nabla E$ accounts for intensive quantities.
Otherwise, the macroscopic energy function can be defined implicitly via Eq.~(\ref{eq:entropy_legendre}) and contact forms~\cite{van2018geometry, van2021classical}.

Remark: the energy function $E$ should be homogeneous of degree 1, so that it verifies for all $\gamma$
\begin{equation}
E(\gamma \,\bm{x}) = \gamma\,E(\bm{x}).
\end{equation}

\subsection{Connection to ports}
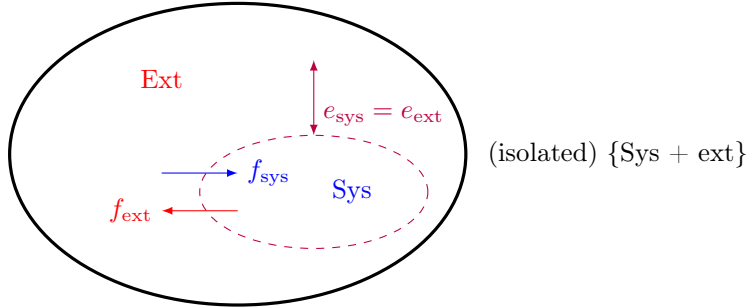
\begin{figure}[h!]
\centering
\begin{tikzpicture}
\draw (0, 0) [very thick]ellipse (3 and 2);
\draw (1, -0.5) [dashed, purple] ellipse (1.5 and 0.75);
\draw (-1, 1) node [red]{Ext};
\draw (1.5, -0.5) node[blue] {Sys};
\draw (5, 0) node {(isolated) \{Sys + ext\}};
\draw (-1, -0.25) [->, >= latex, blue] to ++ (1, 0) [right]node {$f_\mathrm{sys}$};
\draw (-1, -0.75) [left, red] node {$f_\mathrm{ext}$}[<-, >= latex, red] to ++ (1, 0);
\draw (1, 0.25) [<->, >= latex, purple] to ++ (0, 1);
\draw (1, 0.5) [right, purple] node{$e_\mathrm{sys} = e_\mathrm{ext}$};
\end{tikzpicture}
\caption{Flows $f$ and efforts $e$ of a system and its environment. The considered system and its environment as a whole form an isolated system.}
\label{fig:univers}
\end{figure}
The environment acts on the system flow so that at thermodynamic equilibrium, the flows are balanced, and the effort is shared at the system interface (Fig~\ref{fig:univers}):
\begin{equation}
\frac{\partial E^\mathrm{sys}}{\partial\overline{\cF}_i^\mathrm{sys}} = \frac{\partial E^\mathrm{ext}}{\partial\overline{\cF}_i^\mathrm{ext}} \quad \forall i \in \II^x.
\label{eq:effort}
\end{equation}
(see \ref{app:proof_effort} for proof).

Adopting the notations of Eq.~(\ref{eq:SHP}), together with Eq.~(\ref{eq:lambda}), we obtain the relations between flows, efforts and external ports in Table \ref{tab:ports}.

\begin{table}[tb]
\centering
\caption{Port variables and their relations.}
\begin{tabular}{lll}
\toprule
State $\bm{x}$ & $\left[\overline \cS,\, \overline \cF_i\right]$ & $\bm{u} + \dot {\bm{x}} = 0$\\
Effort $\nabla E(\bm{x})$ &$\left[T,\, T\,\lambda_i\right]$ &$\bm{y} = \nabla E(\bm{x})$\\
\bottomrule
\end{tabular}
\label{tab:ports}
\end{table}

\subsection{Conservative, reversible \ac{phs}}

Denoting $\sigma_\mathrm{ext}$ the outgoing entropy flow, the
conservative \ac{phs} interconnection matrix of an open system is found to be
\begin{equation}
\begin{blockarray}{cccccccc}
\begin{block}{cc\BAmulticolumn{3}{c}\BAmulticolumn{3}{c}}
&&\nabla E(\bm{x}) & \bm{u}\\
\end{block}
&&T &\mu &-P &\sigma_\mathrm{ext}  &\dot N_\mathrm{ext} &\dot V_\mathrm{ext}\\
\begin{block}{cc[ccc|ccc]}
        &\dot{\overline{\cS}} &. &.&.&-1 &.&.\\
        \dot{\bm{x}} 
        &\dot {\overline{\cF_n}} &.&. &.&.&-1&. \\
        &\dot {\overline{\cF_v}} &.&.&. &.&.&-1\\
        \BAhline
        &T_\mathrm{ext} &1 &. &.&.&.&.\\
         \bm{y}
         &\mu_\mathrm{ext} &. &1 &. &.&.&.\\
         &-P_\mathrm{ext} &. &. &1 &.&. &.\\
\end{block}
\end{blockarray}.
\label{eq:PHS_thermo}
\end{equation}

\section{Summary of the method and generalization route to irreversible systems}\label{sec:FinalSummary}

This section first summarises the main steps to derive a reversible conservative macroscopic PHS from a microscopic description of matter in an experimental context.
Second,  a process is proposed to complete this modelling under a macroscopic irreversible conservative form, given a dissipation law. 

\subsection{Reversible conservative macroscopic PHS}

In summary, the macroscopic PHS can be described from the microscopic description by completing the following steps:
\begin{enumerate}
\item \textbf{Microstate representation}
Define $\PP$, $\WW = \PP^\star$ and $\MM \subseteq \WW$ equipped with characterizing functions $\mathfrak{F} = \left\{\cF_i:\MM \mapsto \FF_i\right\}_{i\in \II}$.
\item \textbf{Experimental conditions and accessible microstates}
\begin{enumerate}
\item Partition $\mathfrak{F} = \mathfrak{F}^\mathrm{fixed} \cup \mathfrak{F}^\mathrm{free}$ into the set $\mathfrak{F}^\mathrm{fixed}$ of functions the values of which are physically constrained by the experiment and its complement $\mathfrak{F}^\mathrm{free}$, with corresponding sets of indices $\II^\mathrm{fixed}$ and $\II^\mathrm{free}$.
\item Denote $\bm{\theta}^\mathrm{fixed} \subset \text{\Large $\times$}_{i\in \II^\mathrm{fixed}} \FF_i $ the set of experimentally admissible values for functions in $\mathfrak{F}^\mathrm{fixed}$.
\item Denote $\MM_a\left(\bm{\theta}^\mathrm{fixed}\right)$ the corresponding set of admissible microstates.
\end{enumerate}
\item \textbf{Stochastic description}
For all probability distributions $p:\MM_a\left(\bm{\theta}^\mathrm{fixed}\right) \longmapsto [0,\, 1]$, 
\begin{enumerate}
\item Derive the surprisal $\cS^b_p :\bm{m} \in \MM_a\left(\bm{\theta}^\mathrm{fixed}\right) \longmapsto \log_b \frac{1}{p(\bm{m})} \in \RR^+$.
\item Derive the statistical entropy function $\mathsf{S}^k : p \longmapsto \EE_p\left[\cS_p^{b= \exp(\nicefrac{1}{k})}\right] \in \left[0,\, \frac{1}{\Omega}\right]$.
\end{enumerate}
\item \textbf{Boltzmann principle for ergodic systems at thermodynamic equilibrium}
\begin{enumerate}
\item Introduce $\bm{\theta}^\mathrm{free} := \left(\overline{\cF_i}\right)_{i\in \II^\mathrm{free}}$ the values of functions in $\mathfrak{F}^\mathrm{free}$ observed at a macroscopic scale.
\item Define $p^\star\left(\bm{m} \;|\; \mathsf{H}\left(\bm{\theta}^\mathrm{free}\right)\right)$ according to Th.~\ref{th:boltz}.
\item Define the thermodynamic entropy function $S^k :\bm{\theta}^\mathrm{free} \mapsto\mathsf{S}^k\left(p^\star\left(.\;|\; \mathsf{H}\left(\bm{\theta}^\mathrm{free}\right)\right)\right)$.
\end{enumerate}
\end{enumerate}
For common experimental constraints (i.e., constraints on $\mathfrak{F} = \left\{\cF_e,\, \cF_n,\, \cF_v,\, \cS_p^b\right\}$), we obtain the results in Table~\ref{tab:ensembles} (see also \cite{graben1991unified}).

Note that if there is no analytic solution for $\overline \cF_e \mapsto S^k(\overline \cF_e, \dots)$ and its inverse (note that $\dpd{S^k}{\overline \cF_e}$ is monotonic), approximation strategies can be used (see ~\cite{najnudel2021statistical} for an example).

\begin{landscape}
\begin{table*}[tb]
\centering
\captionsetup{width= 10cm}
\caption[Statistical ensembles and associated constraints for usual experimental conditions.]{Statistical ensembles and associated constraints for usual experimental conditions. $\Omega$ denotes the cardinal of $\MMa$ (set of accessible microstates).}
\resizebox{1.5\textwidth}{!}{
\begin{tabular}{p{3cm}|p{3.5cm} p{2.5cm} p{2.5cm} p{5cm} p{4cm} p{5cm}}
\cline{2-7}
&\textbf{Ensemble} &$\bm{\theta}^\mathrm{fixed}$ &$\bm{\theta}^\mathrm{free}$ &$\bm{p^\star}(\bm{m})$ &\textbf{Entropy} &\textbf{Example}\\
\hline
\hline
&Micro-canonical &$\left(E,\,N,\,V,\,S\right)$ & &$\frac{1}{\Omega}$ &$k_B\,\ln\Omega$ &Gas in an isolated tank\\
\cline{2-7}
&\multirow{3}{*}{Isoenthalpic-isobaric} &\multirow{3}{*}{$\left(N,\, S\right)$} &\multirow{3}{*}{$\left(\overline \cF_e,\, \overline \cF_v\right)$} &\multirow{3}{*}{$\frac{1}{\Omega}$} &\multirow{3}{*}{$k_B\,\ln\Omega$} &Gas in a closed tank\\
\multirow{4}{*}{\textbf{No}}&&&&&&with a piston,\\
&&&&&&thermally insulated\\
\cline{2-7}
\multirow{2}{*}{\textbf{thermal contact}}
&\multirow{2}{*}{Unnamed} &\multirow{2}{*}{$\left(V,\, S\right)$} &\multirow{2}{*}{$\left(\overline \cF_e,\, \overline \cF_n\right)$} &\multirow{2}{*}{$\frac{1}{\Omega}$} &\multirow{2}{*}{$k_B\,\ln\Omega$} &Gas in a porous tank,\\
&&&&&&thermally insulated\\
\cline{2-7}
&\multirow{3}{*}{Unnamed} &\multirow{3}{*}{$S$} &\multirow{3}{*}{$\left(\overline \cF_e,\, \overline \cF_n,\, \overline \cF_v \right)$} &\multirow{3}{*}{$\frac{1}{\Omega}$} &\multirow{3}{*}{$k_B\,\ln\Omega$} &Gas in a porous tank\\
&&&&&&with a piston,\\
&&& &&&thermally insulated\\
\hline
\hline
\multirow{10}{*}{\textbf{Thermal contact}}
&\multirow{2}{*}{Canonical} &\multirow{2}{*}{$\left(N,\,V\right)$}&\multirow{2}{*}{$\overline{\cF_e}$} &\multirow{2}{*}{$\frac{\exp\del{-\frac{\cF_e(\bm{m})}{k_B\,T}}}{\cZ(T)}$} &\multirow{2}{*}{$k_B\,\ln \cZ(T)+ \frac{\overline{\cF_e} }{T}$} &Gas in a closed tank,\\
&&&&&& in contact with a thermostat\\
\cline{2-7}
&\multirow{3}{*}{Isothermal-isobaric} &\multirow{3}{*}{$N$} &\multirow{3}{*}{$\left(\overline \cF_e,\,\overline \cF_v\right)$} &\multirow{3}{*}{$\frac{\exp\del{-\frac{\cF_e(\bm{m}) + P\,\cF_v(\bm{m})}{k_B\,T}}}{\cZ(T,\, P)}$} &\multirow{3}{*}{$k_B\,\ln \cZ(T,\, P) + \frac{\overline \cF_e + P\,\overline \cF_v }{T}$} &Gas in a closed tank\\
&&&&&&with a piston,\\
&&&&&&in contact with a thermostat\\
\cline{2-7}
&\multirow{2}{*}{Grand-canonical} &\multirow{2}{*}{$V$} &\multirow{2}{*}{$\left(\overline \cF_e,\,\overline \cF_n\right)$} &\multirow{2}{*}{$\frac{\exp\del{-\frac{\cF_e(\bm{m}) - \mu\,\cF_n(\bm{m})}{k_B\,T}}}{\cZ(T,\, \mu)}$} &\multirow{2}{*}{$k_B\,\ln \cZ(T,\, \mu) + \frac{\overline{\cF_e} - \mu\,\overline\cF_n }{T}$} &Gas in a porous tank,\\
&&&&&& in contact with a thermostat\\
\cline{2-7}
&\multirow{3}{*}{Unnamed} &\multirow{3}{*}{} &\multirow{3}{*}{$\left(\overline \cF_e,\,\overline \cF_n,\,\overline \cF_v\right)$}&\multirow{3}{*}{$\exp\del{-\frac{\cF_e(\bm{m}) + P\,\cF_v(\bm{m}) - \mu\,\cF_n(\bm{m})}{k_B\,T}}$} &\multirow{3}{*}{$\frac{\overline\cF_e + P\, \overline \cF_v -\mu\,\overline \cF_n}{T}$} &Gas in a porous tank\\
&&&&&&with a piston,\\
&&&&&&in contact with a thermostat\\
\bottomrule
\end{tabular}}
\label{tab:ensembles}
\end{table*}
\end{landscape}

\subsection{Irreversible conservative macroscopic PHS from a macroscopic dissipative law}
The system (\ref{eq:PHS_thermo}) models some conservative reversible physics at macroscopic scale.
In some cases, dissipative phenomena can be observed, for which laws are available only at this scale.

In this part, we assume that such a dissipative phenomenon is  described by 
\begin{itemize}
\item[(i)] a  flow-to-effort mapping law 
\begin{equation}
z_{\mathrm{d}}: \bm{f}_{\mathrm{d}} \mapsto \bm{e}_{\mathrm{d}}  = z_{\mathrm{d}}( \bm{f}_{\mathrm{d}} ) 
\text{~~such that for all~}\bm{f}_{\mathrm{d}}, \;\;
z_{\mathrm{d}}(\bm{f}_{\mathrm{d}})^\intercal \bm{f}_{\mathrm{d}}=: P_{\mathrm{d}}\geq 0,
\end{equation}
\item[(ii)] interconnected to the conservative part according to matrix given by
\begin{equation}
\begin{blockarray}{ccccccc}
\begin{block}{cc\BAmulticolumn{2}{c}c\BAmulticolumn{2}{c}}
&&\nabla E(\bm{x}) &z(\bm{w}) &\bm{u}\\
\end{block}
&&T &\bm{e}_\mathrm{s} &\bm{e}_\mathrm{d} &\sigma_\mathrm{ext}  &\bm{e}_\mathrm{ext}\\
\begin{block}{cc[cc|c|cc]}
        \multirow{2}{*}{$\dot {\bm{x}}$} &\dot{\overline{\cS}} &. &. &. &-1 &. \\
        &\bm{f}_\mathrm{s} &. &\bm{J_x} & -\bm{K} &. &-\bm{G_x}\\
        \BAhline
        \bm{w} &\bm{f}_\mathrm{d} &. &\bm{K}^\intercal &\bm{J_w} &. &-\bm{G_w}\\
        \BAhline
        \multirow{2}{*}{$\bm{y}$} &T_\mathrm{ext} &1 &. &. &. &.\\
         &\bm{f}_\mathrm{ext} &. &\bm{G_x}^\intercal &\bm{G_w}^\intercal &. &\bm{J_y}\\      
\end{block}
\end{blockarray},
\label{eq:PHS_diss}
\end{equation}
\item[(iii)] the  dissipated power $P_\mathrm{d}$ being totally converted into an entropy rate 
\begin{equation}
\sigma_{\mathrm{i}} = \nicefrac{P_\mathrm{d}}{T_{\mathrm{d}}}\geq 0 \text{~~where~}T_\mathrm{d}>0,
\end{equation}
denotes the instantaneous macroscopic temperature at which the phenomenon is experienced.
The positivity of $\sigma_{\mathrm{i}}$ reflects the irreversible nature of dissipation.
\end{itemize}

From $z_\mathrm{d}$, we form the irreversible thermodynamic converter with law
\begin{equation}
z: \bm{w} =\left[\bm{f}_{\mathrm{d}}^\intercal,\,T_{\mathrm{d}}\right]^\intercal 
\longmapsto 
\left[ z_{\mathrm{d}}( \bm{f}_{\mathrm{d}} )^\intercal,  
\underbrace{-z_{\mathrm{d}}( \bm{f}_{\mathrm{d}} )^\intercal \bm{f}_{\mathrm{d}}/T_{\mathrm{d}}  
}_{-\sigma_{\mathrm{i}}}
\right]^\intercal 
\end{equation}
where  $-\sigma_{\mathrm{i}}\leq 0$ accounts for the entropy rate incoming into the converter.

This law is conservative as $z(\bm{w})^\intercal \bm{w} = 0$.
Due to irreversibility ($\sigma_{\mathrm{i}}\geq0$), it naturally fulfills the second principle of thermodynamics.

Finally, from (iii), the irreversible conservative thermodynamic macroscopic PHS is given by
\begin{equation}
\begin{blockarray}{cccccccc}
\begin{block}{cc\BAmulticolumn{2}{c}\BAmulticolumn{2}{c}\BAmulticolumn{2}{c}}
&&\nabla E(\bm{x}) &z(\bm{w}) &\bm{u}\\
\end{block}
&&T &\bm{e}_\mathrm{s} &-\sigma_i &\bm{e}_\mathrm{d} &\sigma_\mathrm{ext}  &\bm{e}_\mathrm{ext}\\
\begin{block}{cc[cc|cc|cc]}
        \multirow{2}{*}{$\dot{\bm{x}}$} &\dot{\overline{\cS}} &. &. &-1 &. &-1 &. \\
        &\bm{f}_\mathrm{s} &. &\bm{J_x} &. & -\bm{K} &. &-\bm{G_x}\\
        \BAhline
        \multirow{2}{*}{$\bm{w}$} & T_\mathrm{d} &1 &. &. &. &. &.\\
        &\bm{f}_\mathrm{d} &. &\bm{K}^\intercal &. &\bm{J_w} &. &-\bm{G_w}\\
        \BAhline
        \multirow{2}{*}{$\bm{y}$} &T_\mathrm{ext} &1 &. &. &. &. &.\\
         &\bm{f}_\mathrm{ext} &. &\bm{G_x}^\intercal &. &\bm{G_w}^\intercal &. &\bm{J_y}\\      
\end{block}
\end{blockarray}.
\label{eq:PHS_irr}
\end{equation}

\newpage

\section{Conclusion}
In this paper, we revisited equilibrium \ac{sp} in order to model complex systems with numerous degrees of freedom as macroscopic \ac{phs} with a reduced number of variables.

Starting from the choice of a particle's description and ad hoc characterizing functions, we recalled how to derive the probability of a configuration of particles at equilibrium based on given experimental conditions. In the end, macroscopic variables are revealed to be expectations of the chosen characterizing functions for this probability, and the thermodynamic entropy to be a function of these macroscopic variables. Provided that the energy has been chosen as a characterizing function from the start, the macroscopic energy can in turn be expressed as a function of the thermodynamic entropy and other macroscopic variables. Through the \ac{phs} formalism, experimental conditions are represented as an input flow that acts on the system so that the resulting output is an effort shared with the system. With this formulation, the externality of the environment, as well as its interactions with the system via exchanges of energy and entropy, are made explicit.

As a result, we proposed two \ac{phs} formulations for conservative open systems, a reversible one (with no entropy creation), and an irreversible one (with entropy creation).

An immediate perspective would be to extend this work to non-equilibrium \ac{sp}~\cite{ottinger2005beyond}, so that a macroscopic trajectory would not only be a succession of equilibrium states, and experimental conditions could change faster.

\newpage
\appendix
\section{Proof of Theorem \ref{th:boltz}}
\label{app:proof_boltzmann}
\begin{proof}
\label{proof:poptim}
To solve Eq.~(\ref{eq:entropy_max}), we introduce Lagrange multipliers $\lambda_0$ and $\bm{\lambda}^\mathrm{free} := (\lambda_i)_{i\in \II^\mathrm{free}}$, and optimize~\cite{jaynes1982rationale} the Lagrangian $\mathsf{L}$ defined by
\begin{equation}
\mathsf{L}\; : \del{p,\ \lambda_0,\,  \bm{\lambda}^\mathrm{free}} \longmapsto \mathsf{S}^k(p)
+ \lambda_0\del{\sum_{\bm{m}\in\MM_a\left(\bm{\theta}^\mathrm{fixed}\right)}p(\bm{m}) - 1}
 + \sum_{i \in \II^\mathrm{free}}\lambda_i\del{\EE_p[\cF_i] - \overline{\cF}_i}.
\label{eq:lagrangian}
\end{equation}
A necessary condition to optimize $\mathsf{L}$ is to solve $\frac{\delta\mathsf{L}}{\delta p} = 0$, where $\delta$ denotes the functional derivative.
From Eq.~(\ref{eq:lagrangian})-(\ref{eq:entropy_k})-(\ref{eq:exp}),
\begin{equation*}
\frac{\delta\mathsf{L}}{\delta p} = 0  \Rightarrow \sum_{\bm{m} \in \MM_a\left(\bm{\theta}^\mathrm{fixed}\right)}\del{-k\del{\ln p(\bm{m}) + 1} + \lambda_0 + \sum_{i\in \II^\mathrm{free}}\lambda_i\,\cF_i(\bm{m})} = 0.
\end{equation*}
This is true in particular if $p$ verifies
\begin{equation*}
\begin{split}
&-k\del{\ln p(\bm{m}) + 1} + \lambda_0 + \sum_{i \in \II^\mathrm{free}}\lambda_i\,\cF_i(\bm{m}) = 0 \quad \forall \bm{m} \in \MM_a\left(\bm{\theta}^\mathrm{fixed}\right)\\
\Rightarrow &p(\bm{m}) = \exp\del{\frac{\sum_{i\in \II^\mathrm{free}}\lambda_i\,\cF_i(\bm{m})}{k}}\exp\del{\frac{\lambda_0}{k} - 1} \quad \forall \bm{m} \in \MM_a\left(\bm{\theta}^\mathrm{fixed}\right).
\end{split}
\end{equation*}

A second necessary condition is to solve $\dpd{\mathsf{L}}{\lambda_0} = 0$, which, combined to the first condition, yields
\begin{equation*}
\dpd{\mathsf{L}}{\lambda_0} = 0 \Rightarrow \sum_{\bm{m}\in\MM_a\left(\bm{\theta}^\mathrm{fixed}\right)}p(\bm{m}) = 1\\
\Rightarrow \sum_{\bm{m}\in\MM_a\left(\bm{\theta}^\mathrm{fixed}\right)}\exp\del{\frac{\sum_{i\in \II^\mathrm{free}}\lambda_i\,\cF_i(\bm{m})}{k}} = \exp\del{1 - \frac{\lambda_0}{k}},
\end{equation*}
so that for all $\bm{m} \in \MM_a\left(\bm{\theta}^\mathrm{fixed}\right)$, $p^\star(\bm{m})$ is of the form 
\begin{align}
\widehat{p}^\star\left(\bm{m}\;|\;\mathsf{H}\left(\bm{\theta}^\mathrm{free}\right), \,\bm{\lambda}^\mathrm{free}\right) &= \frac{\exp\del{\frac{\sum_{i\in \II^\mathrm{free}}\lambda_i\cF_i(\bm{m})}{k}}} {\cZ\left(\bm{\lambda}^\mathrm{free}\right)},\\
\text{with } \cZ\left(\bm{\lambda}^\mathrm{free}\right)&:= \sum_{\bm{m}\in\MM_a\left(\bm{\theta}^\mathrm{fixed}\right)}\exp\del{\frac{\sum_{i\in \II^\mathrm{free}}\lambda_i\,\cF_i(\bm{m})}{k}}.
\label{eq:p_hat}
\end{align}

A third necessary necessary condition is to solve $\dpd{\mathsf{L}}{\lambda_i} = 0$ for all $i\in \II^\mathrm{free}$, which, combined with Eq.~(\ref{eq:p_hat}), yields
\begin{equation*}
\begin{split}
\dpd{\mathsf{L}}{\lambda_i} = 0 \Rightarrow &\EE_p[\cF_i] = \overline{\cF}_i\\
\Rightarrow &\frac{\sum_{\bm{m}\in\MM_a\left(\bm{\theta}^\mathrm{fixed}\right)}\cF_i(\bm{m})\exp\del{\frac{\sum_{i\in \II^\mathrm{free}}\lambda_i\,\cF_i(\bm{m})}{k}}}{\cZ\left( \bm{\bm{\lambda}}^\mathrm{free}\right)} = \overline{\cF}_i\\
\Rightarrow &\dpd{}{\lambda_i}k\,\ln\cZ\left( \bm{\bm{\lambda}}^\mathrm{free}\right) = \overline{\cF}_i.
\end{split}
\end{equation*}
We deduce that the optimal distribution $p^\star$ is
\begin{align*}
p^\star\left(\bm{m}\;|\;\mathsf{H}\left(\bm{\theta}^\mathrm{free}\right)\right) &= \widehat{p}^\star\left(\bm{m}\;|\;\mathsf{H}\left(\bm{\theta}^\mathrm{free}\right), \,\bm{\bm{\lambda}}^\mathrm{free}\right),\\
\text{with } \bm{\bm{\lambda}}^\mathrm{free} \text{ such that }\dpd{}{\lambda_i}k\,\ln\cZ\left( \bm{\bm{\lambda}}^\mathrm{free}\right) &= \overline{\cF}_i \quad \forall i \in \II^\mathrm{free}.
\end{align*}
\end{proof}

\section{Proof of effort equality at the system interface}
\label{app:proof_effort}
\begin{proof}
Consider the \emph{isolated} total system constituted by the system under study and its environment.
For all $i \in \II^\mathrm{free}$, we have
\begin{equation*}
\overline{\cF}_i^\mathrm{total} = \overline{\cF}_i^\mathrm{sys} + \overline{\cF}_i^\mathrm{ext}.
\end{equation*}
The entropy is extensive, therefore,
\begin{equation}
S^k_\mathrm{total}\left(\bm{\theta}^\mathrm{free}_\mathrm{total}\right) = S^k_\mathrm{sys}\left(\bm{\theta}^\mathrm{free}_\mathrm{sys}\right) + S^k_\mathrm{ext}\left(\bm{\theta}^\mathrm{free}_\mathrm{ext}\right).
\end{equation}
The total system is isolated, therefore the total entropy is maximal with respect to any variable, so that for all $i \in \II^\mathrm{free}$,
\begin{equation*}
\begin{split}
&\frac{\partial S^k_\mathrm{total}}{\partial\overline{\cF}_i^\mathrm{sys}} = 0\\
\Rightarrow &\frac{\partial S^k_\mathrm{sys}}{\partial\overline{\cF}_i^\mathrm{sys}} + \frac{\partial S^k_\mathrm{ext}}{\partial\overline{\cF}_i^\mathrm{sys}}  = 0\\
\Rightarrow &\frac{\partial S^k_\mathrm{sys}}{\partial\overline{\cF}_i^\mathrm{sys}} - \frac{\partial S^k_\mathrm{ext}}{\partial\overline{\cF}_i^\mathrm{ext}}  = 0\\
\Rightarrow &\lambda_i^\mathrm{sys} - \lambda_i^\mathrm{ext} = 0.
\end{split}
\end{equation*}
Moreover, for all $i \in \II^x$
\begin{equation*}
\frac{\partial E}{\partial\overline{\cF}_i} =  T\,\lambda_i,
\end{equation*}
since
\begin{equation*}
\begin{split}
&S^k\left(E(\overline \cS,\, \overline \cF_i), \, \overline \cF_i \right)
= \overline \cS\\
\Rightarrow &\dpd {S^k}{\overline \cF_e} \dpd {E}{\overline \cF_i} + \dpd {S^k}{\overline \cF_i} = 0\\
\Rightarrow &\frac{1}{T}\dpd {E}{\overline \cF_i} - \lambda_i = 0.
\end{split}
\end{equation*}

We deduce that $\frac{\partial E^\mathrm{sys}}{\partial\overline{\cF}_i^\mathrm{sys}} = \frac{\partial E^\mathrm{ext}}{\partial\overline{\cF}_i^\mathrm{ext}} \quad \forall i \in \II^x$.
\end{proof}

\section*{Acknowledgments}
The authors would like to thank Bernhard Maschke for his comments and helpful discussion.

\bibliography{ref_manuscrit}

\begin{thebibliography}{10}
\expandafter\ifx\csname url\endcsname\relax
  \def\url#1{\texttt{#1}}\fi
\expandafter\ifx\csname urlprefix\endcsname\relax\def\urlprefix{URL }\fi
\expandafter\ifx\csname href\endcsname\relax
  \def\href#1#2{#2} \def\path#1{#1}\fi

\bibitem{boltzmann1877beziehung}
L.~Boltzmann, {\"U}ber die Beziehung zwischen dem zweiten Hauptsatze des
  mechanischen W{\"a}rmetheorie und der Wahrscheinlichkeitsrechnung, respective
  den S{\"a}tzen {\"u}ber das W{\"a}rmegleichgewicht, Kk Hof-und
  Staatsdruckerei, 1877.

\bibitem{landsberg2014thermodynamics}
P.~T. Landsberg, Thermodynamics and statistical mechanics, Courier Corporation,
  2014.

\bibitem{eberard2004port}
D.~Eberard, B.~Maschke, Port-{H}amiltonian systems extended to irreversible
  systems: {T}he example of the heat conduction, IFAC Proceedings Volumes
  37~(13) (2004) 243--248.

\bibitem{eberard2007extension}
D.~Eberard, B.~Maschke, A.~van~der Schaft, An extension of {H}amiltonian
  systems to the thermodynamic phase space: Towards a geometry of nonreversible
  processes, Rep. Math. Phys. 60~(2) (2007) 175--198.

\bibitem{ramirez2013irreversible}
H.~Ramirez, B.~Maschke, D.~Sbarbaro, Irreversible port-{H}amiltonian systems: A
  general formulation of irreversible processes with application to the {CSTR},
  Chem. Eng. Sci. 89 (2013) 223--234.

\bibitem{delvenne2014finite}
J.-C. Delvenne, H.~Sandberg, Finite-time thermodynamics of port-{H}amiltonian
  systems, Physica D 267 (2014) 123--132.

\bibitem{ramirez2016passivity}
H.~Ramirez, Y.~Le~Gorrec, B.~Maschke, F.~Couenne, {On the passivity based
  control of irreversible processes: A port-Hamiltonian approach}, Automatica
  64 (2016) 105--111.

\bibitem{van2021classical}
A.~Van Der~Schaft, Classical thermodynamics revisited: A systems and control
  perspective, IEEE Control Systems Magazine 41~(5) (2021) 32--60.

\bibitem{van2021liouville}
A.~van~der Schaft, Liouville geometry of classical thermodynamics, Journal of
  Geometry and Physics 170 (2021) 104365.

\bibitem{gorban2006model}
A.~N. Gorban, N.~K. Kazantzis, I.~G. Kevrekidis, H.~C. {\"O}ttinger,
  C.~Theodoropoulos, Model reduction and coarse-graining approaches for
  multiscale phenomena, Springer, 2006.

\bibitem{ottinger2007systematic}
H.~C. {\"O}ttinger, {Systematic coarse graining:“Four Lessons and A Caveat”
  from nonequilibrium statistical mechanics}, MRS bulletin 32~(11) (2007)
  936--940.

\bibitem{ising1925beitrag}
E.~Ising, Beitrag zur theorie des ferromagnetismus, Zeitschrift f{\"u}r Physik
  31~(1) (1925) 253--258.

\bibitem{strecka2015brief}
J.~Strecka, M.~Jascur, A brief account of the {I}sing and {I}sing-like models:
  {M}ean-field, effective-field and exact results, arXiv preprint
  arXiv:1511.03031 (2015).

\bibitem{hopcroft2001introduction}
J.~E. Hopcroft, R.~Motwani, J.~D. Ullman, Introduction to automata theory,
  languages, and computation, Acm Sigact News 32~(1) (2001) 60--65.

\bibitem{flajolet2009analytic}
P.~Flajolet, R.~Sedgewick, Analytic combinatorics, Cambridge University Press,
  2009.

\bibitem{liechtenstein1984exchange}
A.~Liechtenstein, M.~Katsnelson, V.~Gubanov, Exchange interactions and
  spin-wave stiffness in ferromagnetic metals, Journal of Physics F: Metal
  Physics 14~(7) (1984) L125.

\bibitem{Grothendieck1986}
A.~Grothendieck, {R\'ecoltes et Semailles (R\'eflexions et t\'emoignage sur un
  pass\'e de math\'ematicien)}, Gallimard (2022), 1986.

\bibitem{plachky2001ideal}
D.~Plachky, An ideal theoretic characterization of finite sets, finite
  algebras, and $\sigma$-algebras of countably generated type, Mathematica
  Slovaca 51~(3) (2001) 301--311.

\bibitem{gray2011entropy}
R.~M. Gray, Entropy and information theory, Springer, 2011.

\bibitem{knuth1985dynamic}
D.~E. Knuth, Dynamic {H}uffman coding, Journal of algorithms 6~(2) (1985)
  163--180.

\bibitem{patrascioiu1987ergodic}
A.~Patrascioiu, The ergodic-hypothesis: {A} complicated problem in mathematics
  and physics, Los Alamos Science 15 (1987) 263--279.

\bibitem{zia2009making}
R.~K. Zia, E.~F. Redish, S.~R. McKay, {Making sense of the Legendre transform},
  American Journal of Physics 77~(7) (2009) 614--622.

\bibitem{van2018geometry}
A.~Van~der Schaft, B.~Maschke, Geometry of thermodynamic processes, Entropy
  20~(12) (2018) 925.

\bibitem{ottinger2005beyond}
H.~C. {\"O}ttinger, {Beyond Equilibrium Thermodynamics}, John Wiley \& Sons,
  2005.

\bibitem{duindam2009modeling}
V.~Duindam, A.~Macchelli, S.~Stramigioli, H.~Bruyninckx, Modeling and control
  of complex physical systems: {T}he port-{H}amiltonian approach, Springer
  Science \& Business Media, 2009.

\bibitem{van2014port}
A.~J. van~der Schaft, D.~Jeltsema, et~al., Port-{H}amiltonian systems theory:
  {A}n introductory overview, Foundations and Trends in Systems and Control
  1~(2-3) (2014) 173--378.

\bibitem{falaize2016passive}
A.~Falaize, T.~H{\'e}lie, Passive guaranteed simulation of analog audio
  circuits: {A} port-{H}amiltonian approach, Applied Sciences 6~(10) (2016)
  273.

\bibitem{remy2021time}
R.~M{\"u}ller, Time-continuous power-balanced simulation of nonlinear audio
  circuits: Realtime processing framework and aliasing rejection, Ph.D. thesis,
  Sorbonne Universit{\'e} (2021).

\bibitem{bigelow2020power}
T.~A. Bigelow, Power and energy in electric circuits, in: Electric Circuits,
  Systems, and Motors, Springer, 2020, pp. 105--121.

\bibitem{graben1991unified}
H.~Graben, J.~R. Ray, Unified treatment of adiabatic ensembles, Phys. Rev. A
  43~(8) (1991) 4100.

\bibitem{najnudel2021statistical}
J.~Najnudel, T.~H{\'e}lie, D.~Roze, R.~M{\"u}ller, From statistical physics to
  macroscopic port-{H}amiltonian systems: {A} roadmap, in: 7th IFAC Workshop on
  Lagrangian and {H}amiltonian Methods for Nonlinear Control, 2021.

\bibitem{jaynes1982rationale}
E.~T. Jaynes, On the rationale of maximum-entropy methods, Proceedings of the
  IEEE 70~(9) (1982) 939--952.

\end{thebibliography}

\end{document}